\def\N{I\!\!N}

\def\Proof{{\bf Proof.~}}

\def\|{\,\parallel\,}

\def\Theoreme#1{
\begin{tabular}{p{0.1cm}|p{13.7cm}}&
\vspace{0cm} 
\begin{#1}
}
\def\endTheoreme#1{
\end{#1}
\end{tabular} 
\\
}

\newtheorem{theo}{Theorem}
\def\endth{\end{theo}}
\newtheorem{defi}{Definition}
\newcommand\defn{\begin{defi}}
\def\enddefn{\end{defi}}

\newtheorem{propo}{Proposition}
\newcommand\prop{\begin{propo}}
\def\endprop{\end{propo}}

\newtheorem{coro}{Corollary}
\newcommand\cor{\begin{coro}}
\def\endcor{\end{coro}}

\newtheorem{lemmes}{Lemma}
\newcommand\lemme{\begin{lemmes}}
\def\endlemme{\end{lemmes}}
\newtheorem{demonstration}{Demonstration}
\newcommand\dem{\begin{demonstration}}
\def\enddem{\end{demonstration}}
\newtheorem{propriete}{Propri\'{e}t\'{e}}
\newcommand\propri{\begin{propriete}}
\def\endpropri{\end{propriete}}

\documentstyle [epsf, 11pt]{article}

\title{The Structure of Chip Firing Games and related Models
\thanks{This work have been done during the time of M.M and H.D.P in 
Departamento de Ingenier\'{\i}a 
Matem\'atica, Universidad de Chile and was supported by Project ECOS-C96E02
and Chilean program FONDAP in Applied Mathematics (E.G., M.M., H.D.P.)
}}
\author{Eric Goles\thanks{Departamento de Ingenier\'{\i}a Matem\'atica,
Escuela de Ingenier\'{\i}a, Universidad de Chile, Casilla 170-Correo 3,
 Santiago, Chile.
Email : egoles@dim.uchile.cl} \quad
  Michel Morvan\thanks{LIAFA Universit\'{e} Denis Diderot 
Paris 7 and Institut universitaire de France - Case 7014-2, Place Jussieu-75256 Paris Cedex 05-France. 
Email: morvan@liafa.jussieu.fr}\quad
 Ha Duong Phan\thanks{LIAFA Universit\'{e} Denis Diderot 
Paris 7 - Case 7014-2, Place Jussieu-75256 Paris Cedex 05-France. 
Email: phan@liafa.jussieu.fr} \quad 
}

\date{} 

\begin{document}

\maketitle 

\newpage

{\bf Abstract:} {\em In this paper, we study the dynamics of sand
grains falling in sand piles. Usually sand piles are characterized by
a decreasing integer partition and grain moves are described in terms
of transitions between such partitions. We study here four main
transition rules.  The more classical one, introduced by Brylawski
\cite{Bry73} induces a lattice structure $L_B (n)$ (called dominance
ordering) between decreasing partitions of a given integer $n$. We
prove that a more restrictive transition rule, called $SPM$ rule,
induces a natural partition of $L_B (n)$ in suborders, each one
associated to a fixed point for $SPM$ rule. In the second part, we
extend the $SPM$ rule in a natural way and obtain a model called Chip
Firing Game \cite{GK93}. We prove that this new model has interesting
properties: the induced order is a lattice, a natural greedoid can be
associated to the model and it also defines a strongly convergent
game. In the last section, we generalize the $SPM$ rule in another way
and obtain other lattice structure parametrized by some $\theta$:
$L(n,\theta)$, which form for $\theta \in [-n+2,n]$ a decreasing
sequence of lattices. For each $\theta$, we characterize the fixed
point of $L(n,\theta)$ and give the value of its maximal sized chain's
lenght. We also note that $L(n,-n+2)$ is the lattice of all
compositions of $n$.  }
 
\newpage

\section{Introduction}
The set of all partitions (decreasing compositions) of a given integer $n$ has been extensively investigated \cite{And76}. The orders defined on it play an important role, especially the  so-called dominance ordering, following the result of Brylawski \cite{Bry73}, which shows the lattice structure of the object (denoted by $L_B (n)$). The dominance ordering $\leq_B$ on partitions of $n$ is defined by: if $a=(a_1,\ldots,a_n)$ and $b=(b_1,\ldots,b_n)$ are decreasing compositions of $n$, i.e. $a_1 \geq a_2 \geq \ldots \geq a_n$, $b_1 \geq b_2 \geq \ldots \geq b_n$ and $a_1 + \ldots + a_n = b_1 + \ldots + b_n = n$, then $a\leq_B b$ if $\forall i, 1 \leq i \leq n, \sum_{j=1}^{j=i} a_j \leq  \sum_{j=1}^{j=i} b_j$. Since then, some results on the structure of this lattice have been presented, e.g. those concerned with maximal chains, fixed points or about its general  structure \cite{Bry73,GK93,GK86}.

An interesting application of this problem is the sand piles' problem, which has been investigated in many works in other physics and combinatorics \cite{ALSSTW89,Bak88,Bry73,GK93,Spe86}. The core of this problem is to study a model of the sand piles, corresponding to the partitions of a certain integer, and the possible moves to transform one sand pile to another.
In this context, two main models have been investigated until now. The first one is the model in which two falling rules allow to obtain any decreasing partition of $n$ starting from the partition $N = (n,0,\ldots,0)$. The obtained order is exactly $L_B (n)$, which is also the name of the model \cite{Bry73}. In the second one, called $SPM(n)$ \cite{GK93,GMP97} for Sand Piles Model, only one rule is kept, which induces a suborder $\leq_{SPM}$ of $L_B (n)$, that is also a lattice, denoted by $SPM(n)$. The two falling rules discribing $L_B (n)$ are called the horizontal rule and the vertical rule, applicable to a decreasing sequence $a$, $a=(a_1,\ldots, a_n)$ such that $\sum a_i =n$, (see Figure 1 and Figure 2).
\begin{trivlist}
\item[] {\em Vertical rule:} $a_{1},\ldots,a_{i},a_{i+1},\ldots, a_{n}
\rightarrow a_{1},\ldots,a_{i}-1,a_{i+1}+1,\ldots,a_{n} $; \\
this rule can be applied when $a_{i} - a_{i+1} \geq 2$. 
\item[] {\em Horizontal rule:} $a_1,\ldots, p+1,p,\ldots,p,p-1,\ldots,a_n
\rightarrow a_1,\ldots,p,p,\ldots,p,p,\ldots,a_n$.
\end{trivlist}
 $SPM(n)$ is obtained by starting from the partition $(n,0,\ldots,0)$ and applying only the vertical rule (also called $SPM$ rule).

In Section 2, we investigate the problem of the structure of the set of fixed points obtained when we repeatedly apply the vertical rule starting from any element of $L_B (n)$. We show that this set is a lattice  anti-isomophic to the lattice of strictly decreasing partitions of $n$ ordered by the dominance ordering. Moreover, we show that these fixed points induce a natural partition of $L_B (n)$. 

In the next section, we investigate a natural way of extending $SPM$ by modifying the vertical rule in the following way. For a fixed integer $m$, at each step, we allow $m$ (instead of one) sand grains to fall at the same time, each one arriving on one of the $m$ columns on the right. This rule can only be applied if the obtained partition is still decreasing. We prove that this new model, called the Chip Firing Game, has interesting properties: the induced order is a lattice, a natural greedoid can be associated to it and it also defines a strongly convergent game. These results show that these kinds of objects, which originally come from physical considerations, can be interesting in many different ways: order and lattice theory, game theory, language theory.

Another extension of the $SPM(n)$  model based on an extension of the conditions of possible moves is studied in Section 4. In the original $SPM(n)$ model, in order to have a move, the difference in height of two consecutive piles must be greater or equal to 2. If, instead of 2, we consider other values, denoted by $\theta$, we will obtain new models that we will denote by $L(n,\theta)$. In Section 4, we prove that for any $\theta$, $L(n,\theta)$ is a lattice and furthermore, that the order of $L(n,\theta)$ corresponds to the dominance ordering. We also show that for any $\theta \in [-n+2,n-1]$, $L(n,\theta +1)$ is a suborder of $L(n,\theta)$. Moreover, the lattice $L_B (n)$ is shown to be a sublattice of $L(n,1)$. At the end of the section, we present an explicit formula for fixed points and maximal sized chains' lenght for each lattice $L(n,\theta)$.

\bigskip In the following, we are going to discuss about some lattice properties of the above dynamical systems. Let us recall that a finite lattice can be described as a finite partial order such that any two elements $a$ and $b$ admit a least upper bound (denoted by $sup(a,b)$) and a greatest lower bound (denoted by $inf(a,b)$). $Sup(a,b)$ is the smallest element among the elements greater than both $a$ and $b$. $Inf(a,b)$ is defined similarly. A useful result about  finite lattices is that a partial order is a lattice if and only if it admits a greateast element and any two elements admit a greatest lower bound. For more informations about lattice theory, see \cite{Bir67,Dav90}.

\begin{figure}[htb]
$${\epsfbox{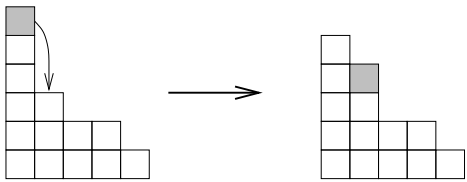}}$$
\caption{The movement of sand grain by vertical and horizontal rules}
\label{fig:01}
\end{figure} 

\begin{figure}[htb]
$${\epsfbox{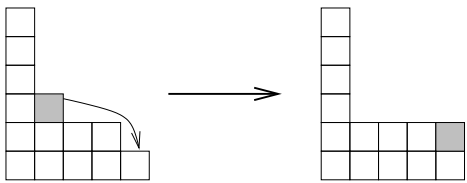}}$$
\caption{The movement of sand grain by horizontal rule}
\label{fig:02}
\end{figure}
%
\section{Characterization of all fixed points by $SPM$ moves and a partition of $L_B (n)$}
In this part we are going to show that the set of all elements of $L_B (n)$ on which $SPM$ rule cannot be applied (called {\em $SPM(n)$-fixed points}) is in bijection with the set of all strictly decreasing partitions (called {\em strict partitions} in the following) of $n$ and that the two induced lattices are anti-isomorphic. Moreover, we will show that these fixed points induce a natural partition of $L_B (n)$.

Let us first denote by $\Phi$ the set of all $SPM(n)$-fixed points and by $StrictPar$ the set of all strict partitions of $n$. Let us recall the notion of dual or conjugate application \cite{Bry73}: for a given partition $a=(a_1,\ldots,a_n)$, the {\em dual} $d(a)$ is the partition $d(a) =  (a^*_1,\ldots a^*_n)$ such that $a^*_i = |\{a_j, a_j \geq i\}|$. $d(a)$ is also called the {\em conjugate} of $a$ and denoted by $a^*$. 
\begin{theo}
The set $\Phi$ of all $SPM(n)$-fixed points is in bijection with the set $StrictPar$ of all strict partitions by the dual application $d$. Moreover $\Phi$ is a lattice  with the maximal element being the fixed point $P_0$ of $SPM(n)$ (starting from $N=(n,0,\ldots,0)$), the minimal element being $1=(1,\ldots,1)$ and the dual application $d$ is an anti-isomorphism.
\end{theo}
\Proof 
Let us first consider a property of the dual application: given a partition $a=(a_1,\ldots,a_n)$, by definition $a^*_i = |\{a_j, a_j \geq i\}|$, so we have $|\{a_j, a_j =i\}| = a^*_{i}-a^*_{i+1}$. A fixed point $P=(p_1,\ldots,p_n)$ by $SPM$ rule is a partition such that $\forall j$ $ p_j - p_{j+1} \leq 1$, which implies that $\forall i \leq p_1$ $\{p_j, p_j =i\} \neq \emptyset$, and that $p^*_{i}-p^*_{i+1} \geq 1$, so $d(P)$ is a strict partition. By the same way, we can prove that the dual of a strict partition is a fixed point by $SPM$ rule, and the dual application is then a bijection between $\Phi$ and $StrictPar$.

On the other hand, let us recall a remark of Greene and Kleitman \cite{GK86}: $a \rightarrow b$ is a vertical movement (by $SPM$ rule) if and only if  $d(b) \rightarrow d(a)$ is a horizontal movement, which implies that $d$ is a anti-isomorphism, i.e. $P_1\geq_B  P_2 \Leftrightarrow d(P_1) \leq_B d(P_2)$. Moreover, the smallest strict partition (by dominance ordering) of $n$ if $n=\frac{1}{2}k(k+1)+k', 0 < k' \leq k$, is clearly equal to $(k+1,k,\ldots,k+2-k',k-k',k-k'-1,\ldots,1,0,\ldots,0)$ which corresponds to the fixed point $P_0= (k,k-1,\ldots,k'+1,k',k',k'-1,\ldots 1,0,\ldots,0)$, which is the fixed point of $SPM(n)$. In order to prove that $\Phi$ is a lattice we will prove that   $StricPar$ is a lattice. For that, it is sufficient to show that for two given strict partitions $a$ and $b$, the partition $c=inf_{L_B (n)}(a,b)$ is also a strict partition. Let us suppose that $c$ is not a strict partition, i.e. $ \exists i$ $c_i=c_{i+1}$, knowing that $\forall m $ $\sum_{j=1}^{j=m} c_j = min (\sum_{j=1}^{j=m} a_j, \sum_{j=1}^{j=m} b_j)$, and without loss of generality, we may suppose that $\sum_{j=1}^{j=i-1} c_j= \sum_{j=1}^{j=i-1} a_j$. Two cases are now possible: either $\sum_{j=1}^{j=i} c_j= \sum_{j=1}^{j=i} a_j$, or $\sum_{j=1}^{j=i} c_j= \sum_{j=1}^{j=i} b_j$. In the first case, we obtain $a_i = c_i$, and as  $a$ is a strict partition, then $a_{i+1} < a_i = c_i = c_{i+1}$, and then $\sum_{j=1}^{j=i+1}a_j < \sum_{j=1}^{j=i+1} c_j$, which is a contradiction. The second case is similar, so $c$ is a strict partition, which ends our proof.
\hfill$\Diamond$\\ 
Using the previous results on the characterizations of fixed points, we can find a partition of $L_B (n)$ by orders that satisfy SPM(n) rule. First consider the set $U_{SPM(n)}(P)$ containing all elements $a$ of $L_B (n)$ such that $a\geq_{SPM}  P$. For every two different $P_1$, $P_2$ of $\Phi$, $U_{SPM(n)}(P_1) \cap U_{SPM(n)}(P_2) = \emptyset$, i.e.
\begin{lemmes} For every pair of fixed points $P_1$, $P_2$, there does not exist any element $a \in L_B (n)$ such that $a \geq_{SPM} P_1$ and $a \geq_{SPM} P_2$.
\end{lemmes}
\Proof:
Suppose that such an element $a$ exists. Consider the interval $\{c \in L_B(n) | a \geq_B c \geq_B (1,\ldots,1)\}$. From Greene and Kleitman \cite{GK86}, we know that there exists only one $b$ in $L_B (n) $ such that $b$ is the minimal in $L_B (n)$ satisfying $a\geq_{SPM} b\geq_{B} 1$. Since $P_1$ satisfies this condition, $P_1=b$. Moreover, since such an element $b$ is unique, $P_2=b$, which is a contradiction. 
\hfill$\Diamond$\\
Now, on the other hand, for each element $a$ of $L_B (n)$ there always exist a $P$ of $\Phi$ such that $a\geq_{SPM} P$, so we obtain the following result, illustrated on Figure \ref{fig:03}: 
\begin{figure}
$${ \epsfbox{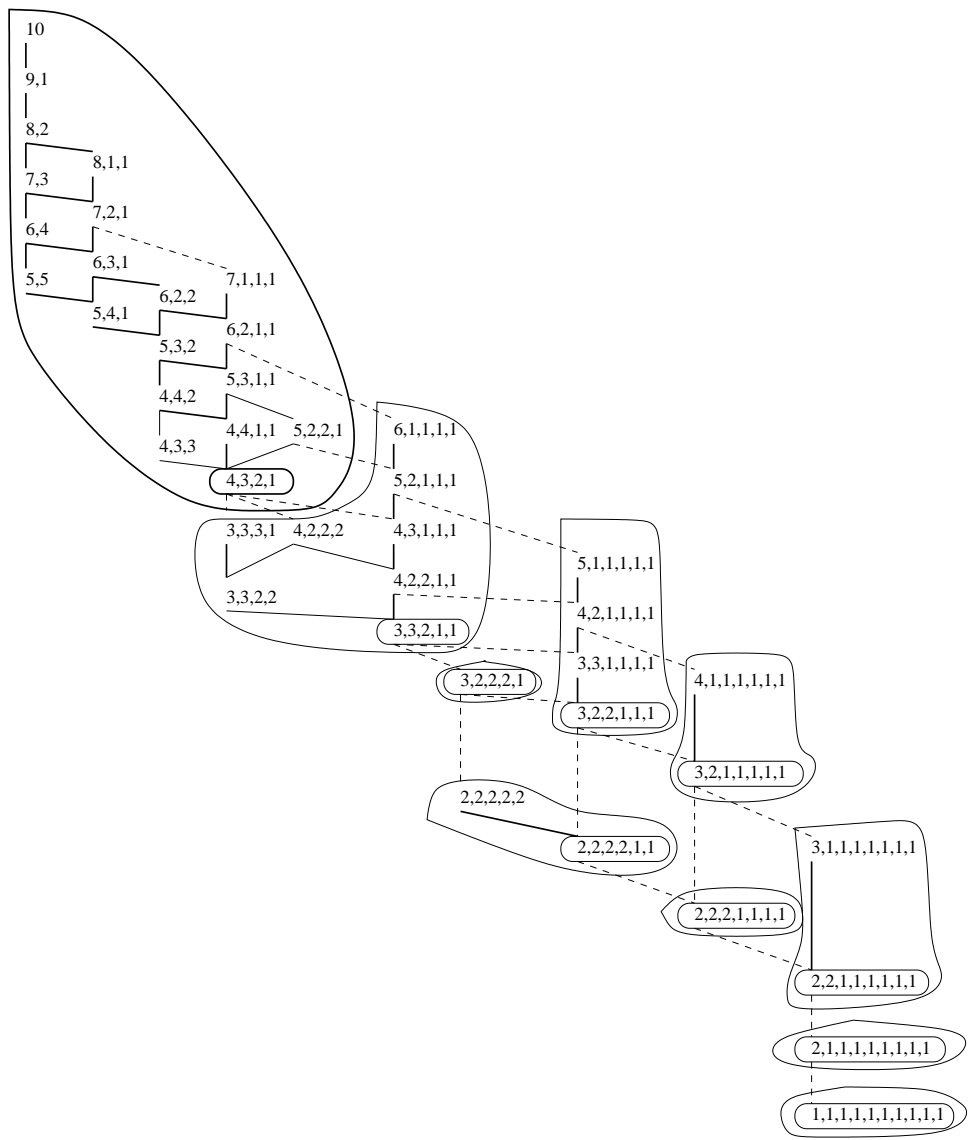}}$$
\caption{The partition of $L_B (n)$ by $\{U_{SPM(n)}(P),P \in \Phi\}$ in the case $n=10$: all fixed points are surrounded, and full lines signify vertical movements (transitions by $SPM$ rule), while dashed lines signify horizontal movements.}
\label{fig:03}
\end{figure}
\begin{theo} $\{ U_{SPM(n)}(P), P \in \Phi \}$ is a partition of $L_B (n)$, i.e.
$L_B (n) = \bigcup U_{SPM(n)}(P)$ and for all $P_1 \neq P_2$, $U_{SPM(n)}(P_1) \cap U_{SPM(n)}(P_2) = \emptyset.$
\end{theo}

\section{The structure of $CFG(n,m)$}
In this section, we define new rules for grains to move, obtaining the model $CFG(n,m)$ and we give some results on the structure of this model. Some notions defined in this section have been first introduced in \cite{Eri93}. The section is organized as follows: we first define the falling rules, then characterize the condition for a configuration to be obtainable from another, we show that the order naturally associated to $CFG(n,m)$ is a lattice and we prove that $CFG(n,m)$ can also be seen as a strongly convergent game or, in the language theory context, as a greedoid.

Let $n$ and $m$ be two integers, $CFG(n,m)$ is a chip firing game containing partitions of $n$ with the update rule defined as follows \cite{GK93} (see Figure 4):
$$ \left\{\begin{array}{ll}
	  x_i \rightarrow x_i - m  \\
	  x_j \rightarrow x_j + 1 & \forall j \in \{i+1,\ldots i+m \}
	\end{array}
\right. $$
with the condition: $x_i - x_{i+1} \geq m+1$.
\begin{figure}
$${\epsfbox{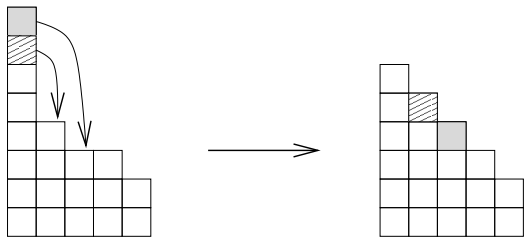}}$$
\caption{The movement of sand grain by rule CFG(n,2)}
\label{fig:04}
\end{figure} 

Let $a$ and $b$ be two partitions of $n$, we say that $b \leq_{(n,m)} a$ if we can apply a sequence of admissible transitions from $a$ to obtain $b$. Let $O$ be a given partition. We define the order $CFG(O,n,m)$ as the suborder of $CFG(n,m)$ induced by all elements of $CFG(n,m)$ smaller than $O$ for $\leq_{(n,m)}$. Remark that  $CFG(N,n,m)$ = $CFG(n,m)$.

The first result of this part is analogous to the one for the  other classes of models:
\begin{theo}
$CFG(O,n,m)$ is a lattice.
\end{theo}
In order to prove  this theorem we need to introduce some new objects and to prove intermediate results. 
Let us first consider a sequence of admissible transitions from $O$ to $a$:
$$O=(O_1,O_2,\ldots,O_n) \rightarrow   \ldots \rightarrow a=(a_1,a_2,\ldots,a_n)$$
The {\em shot vector} (see for example \cite{Eri93}) of this sequence is the vector $k(O,a) \in \N^n$ whose entry $k_i$ is  the number of times we applied the rule to column $i$  during the  sequence. In $CFG(O,n,m)$, it is easy to see that $k(O,a)$ is defined by the following formula, regadless of the choice of transition sequence:
$$k_i = \frac{O_i - a_i + k_{i-m} + \ldots + k_{i-1}}{m}$$
in which $k_j$ is defined to be $0$ if $ j \leq 0$.

Let us denote by $|k(O,a)|$ the sum $\sum_{i=1}^{i=n} k_i (O,a)$.
Let $a$ and $b$ be two partitions of $CFG(O,n,m)$, we define that $k(O,a) \leq k(O,b)$ if $\forall i$ $ k_i(O,a) \leq k_i(O,b)$. Moreover, if $a \geq_{(n,m)} b$, it is clear that $k(O,b) = k(O,a)+k(a,b)$. Let us give here a useful result about the shot vector:
\begin{lemmes}
Let $a$ and $b$ be two partitions of $CFG(O,n,m)$ such that there exists an index $j$ such that
 $k_{j}(O,a) \leq k_{j}(O,b)$ and $ \forall j^{'} \neq j , k_{j^{'}}(O,a) \geq k_{j^{'}}(O,b)$. If it is possible to apply the transition $CFG(n,m)$ at the position $j$ of $b$, then it is also possible to apply this transition to $a$ at the same position.
\end{lemmes}
\Proof
Knowing that the necessary and sufficient condition to apply the transition at the position $j$ of $b$ is $ b_{j} - b_{j+1} \geq m+1 $, let us consider the difference $a_{j} - a_{j+1} $. As 
$$ k_{j}(O,a) = \frac{ O_{j} - a_{j} + k_{j-m}(O,a)+ \ldots + k_{j-1}(O,a) } {m}, $$ 
$$ a_{j} - a_{j+1} = (O_{j}-O_{j+1}) + k_{j-m}(O,a) - (m+1)k_{j}(O,a) + mk_{j+1}(O,
a)$$
so $$(a_{j} - a_{j+1})-(b_{j} - b_{j+1}) = $$
$$ (k_{j-m}(O,a) - k_{j-m}(O,b)) + 
m(k_{j+1}(O,a) - k_{j+1}(O,b)) - (m+1)(k_{j}(O,a) - k_{j}(O,b))$$ $$ \geq 0, $$
  which implies that  $$a_{j} - a_{j+1} \geq m+1,$$
which proves the result.
\hfill$\Diamond$\\
We can now characterize the order relation between elements of $CFG(O,n,m)$ (Proposition 1) and the fomula for $inf_{(n,m)}(a,b)$ for two given elements $a$ and $b$ of $CFG(O,n,m)$ (Proposition 2).
\begin{propo}
If $a$ and $b$ are two partitions of $CFG(O,n,m)$, then:
$$a \geq_{(n,m)} b \Longleftrightarrow k(O,a) \leq k(O,b).$$
\end{propo}
\Proof
If $a \geq_{(n,m)} b$ then $k(O,b) = k(O,a)+k(a,b) \geq k(O,a)$. Let us now assume that $k(O,a) \leq k(O,b)$  and consider two sequences of $CFG(O,n,m)$ transitions, one from $O$ to $a$ and other from $O$ to $b$:   
$$ O \rightarrow c_{1} \rightarrow \ldots \rightarrow c_{r} \rightarrow a$$
$$ O \rightarrow d_{1} \rightarrow \ldots \rightarrow d_{s} \rightarrow b.$$
We will construct step by step a sequence of transitions: $ a \rightarrow e_{1} \rightarrow \ldots \rightarrow e_{t} \rightarrow b$ showing that $a \leq_{(n,m)} b$. Knowing that $(0,0,\ldots,0) =  k(O,O) \leq k(O,a) \leq k(O,b) $, there exists a first index $i$ such that there exists $j$ such that $ k_{j}(O,d_{i}) > k_{j}(O,a)$ and $ \forall j^{'} \neq j , k_{j^{'}}(O,d_{i}) \leq k_{j^{'}}(O,a)$. Since $k(O,d_{i-1}) \leq k(O,a)$ then $k_{j}(O,d_{i-1}) = k_{j}(O,a)$ and $ k_{j}(O,d_{i}) = k_{j}(O,a) + 1 $. Since $d_{i-1}$ and $a$ satisfy the conditions of Lemma 2, we can apply the transition $CFG$ at the position $j$ of $a$ to obtain a new partition, denoted by $e_{1}$, and we have $ k(O,d_{i}) \leq k(O,e_{1}) \leq k(O,b)$. By repeating this procedure, we can define  $e_{2}, e_{3}, \ldots $ Since $ | k(O,e_{l}) - k(O,a) | = l$ and $ k(O,a) \leq k(O,e_{l}) \leq k(O,b)$ then, after $ t = |k(O,b) - k(O,a) |$ steps, we will have $ e_{t} \rightarrow b$. A sequence of transitions $CFG(n,m)$ from $a$ to $b$ is then established.
\hfill$\Diamond$\\
\begin{propo}
Let $a$ and $b$ be two configurations of $CFG(O,n,m)$. Let $k=(k_1,k_2,\ldots,k_n)$ such that for each $i$, 
$k_i= max(k_i (O,a), k_i (O,b))$.
Then the configuration $c$  such that $k(O,c) =k$ is in $CFG(O,n,m)$ and $c = inf_{(n,m)}(a,b)$.
\end{propo}
\Proof
In order to prove that $c=inf_{(n,m)}(a,b)$, we are going to show that $a \geq_{(n,m)} c$ and $b \geq_{(n,m)} c$. Since $c$ is cleary the greatest partition that can satisfy these properties, this will prove the result. Let us assume that $k(O,a)$ and $k(O,b)$ are not comparable (otherwise, $a$ and $b$ are comparable and the result is obvious). Let us show that $a \geq_{(n,m)} c$ (the proof is similar for $b$). For that,  it is sufficient to find a partition $a'$ such that $a \rightarrow a'$ and $k(O,a') \leq k(O,c)$. We are going to prove the existence of such a partition by using a sequence from $O$ to $b$. Let  $O \rightarrow d_{1} \rightarrow \ldots \rightarrow d_{s} \rightarrow b$ be such a sequence and let $l$ be the first index such that $k(O,d_l) \leq k(O,a)$ and $k(O,d_{l+1}) \not\leq k(O,a)$. Let us consider the position $i$ at which the transition is applied for $d_l$. We have $k_i(O,d_l) \leq k_i(O,a)$ and $k_i(O,d_{l+1}) > k_i(O,a)$. Since $a$ and $d_l$ satisfy the conditions of Lemma 2,  we can apply the transition at position $i$ of $a$ to obtain a new partition $a'$. The shot vector of $a'$ satisfies $\forall j \neq i$ $k_j(O,a') = k_j(O,a) \leq k_j(O,c)$ and $k_i(O,a') = k _i(O,d_{l+1}) \leq k_i(O,b) \leq k_i(O,c)$, which were the wanted conditions. This proves the result.
\hfill$\Diamond$\\
The well-known fact that an order is a lattice if it contains a maximal element and if it is closed by the $inf$ (see for example \cite{Dav90}) give us immediately Theorem 3.

As we have said in the introduction of this section, $CFG(n,m)$ presents other interesting facets. Let us now consider $CFG(n,m)$ as a game where the play rule is the transition rule of $CFG(n,m)$. By means of the result of Proposition 2, we will prove that $CFG(n,m)$ is a strongly convergent game. Let us now  give the definition of such games:
\begin{defi}\cite{Eri93}
A game is said to have the strongly convergent property if, given any starting position, either every play sequence can be continued indefinitely, or every play sequence will converge to the same terminal position in the same number of moves.
\end{defi}
So we have the following result, which is a corollary of Theorem 3.
\begin{coro}
$CFG(n,m)$ is a strongly convergent game.
\end{coro}
\Proof
Let $O$ be a partition and let $O \rightarrow a_1 \rightarrow a_2 \rightarrow \ldots $ and $O \rightarrow b_1 \rightarrow b_2 \rightarrow \ldots $ be two play sequences from $O$. From Theorem 3, $CFG(O,n,m)$ is a lattice, so the two sequences converge to the same terminal configuration, denoted by $c$. Moreover, the length of a sequence from $O$ to $c$ is equal to $\sum k_{i}(O,c)$, so all play sequences from $O$ to $c$ have the same length and $CFG(n,m)$ is then a strongly convergent game. 
\hfill$\Diamond$\\
Let us now see that the chip firing game can induce interesting properties in terms of  language. $CFG(O,n,m)$  defines  a language in the following way: the alphabet of $CFG(O,n,m)$ is $A = \{1,2,\ldots ,n\}$, and a word $\alpha=(i_1,i_2,\ldots, i_r) \in A^{*}$ belongs to $CFG(O,n,m)$ if there is  a play sequence $O = c_0 \rightarrow c_1 \rightarrow \ldots \rightarrow c_r =a$ such that $c_{s-1} \rightarrow c_s$ is a transition at position $i_s$. 

\begin{defi} \cite{KLS91}
A language $L$ is a {\em greedoid} if it is left-hereditary which means that:
$ \alpha \gamma \in L \Rightarrow \alpha \in L$
and if $L$ also satisfies the following exchange condition:
$$\alpha,\beta \in L, |\beta| > |\alpha| \Rightarrow \mbox{ there exists a letter $x$ in $\beta$ such that } \alpha x \in L.$$
\end{defi}
We have:
\begin{theo}
For all $O,n,m$, $CFG(O,n,m)$ is a greedoid.
\end{theo} 
\Proof
Since it is clear that $CFG(O,n,m)$ satisfies the left-hereditary condition, let us now prove that it satisfies the exchange condition. Let $\alpha$ and $\beta$ be two words of $CFG(O,n,m)$ such that  $|\beta| > |\alpha|$. Let $a$ and $b$ be the corresponding partitions of $\alpha$ and $\beta$ respectively.  Using the same argument as in Proposition 2, we can find a new word $\alpha '$ and its corresponding partition  $a'$ such that $a \rightarrow a'$ and $\alpha ' = (\alpha,i)$ where $i$ is a letter in $\beta$. The theorem is then proved.
\hfill$\Diamond$\\
\section{Another extension allowing non decreasing compositions of $n$}
In \cite{GMP97}, we studied the models $SPM(n)$ and $L_B (n)$ and the structure of the set of all partitions for a given integer $n$. Futhermore, in the previous sections, we studied an extension model of $SPM(n)$, which also contains  partitions.

 Let us now consider the set of all compositions of $n$, that is, a sequence of $n$ integers, the sum of which is equal to $n$, i.e. $(a_1, \ldots, a_n), a_i \geq 0, \sum a_i = n$. In $SPM(n)$, each time we apply the transition $a_{1},\ldots,a_{i},a_{i+1},\ldots, a_{n} \rightarrow a_{1},\ldots,a_{i}-1,a_{i+1}+1,\ldots,a_{n} $ we need the condition that the difference $a_i -a_{i+1}$ is at least $2$ which implies that the new sequence is still decreasing. In this section, without the condition of decreasing sequence, instead of the value 2, we can consider other values, denoted by $\theta$, to obtain different new sets of compositions. Let us denote by $T(\theta)$ a transition $a_{1},\ldots,a_{i},a_{i+1},\ldots, a_{n} \rightarrow a_{1},\ldots,a_{i}-1,a_{i+1}+1,\ldots,a_{n}$ with the condition that $a_i - a_{i+ 1} \geq \theta$ and by $L(n,\theta)$ the set of all compositions that we can obtain by iterating rule $T(\theta)$, starting from $N=(n,0,\ldots,0)$. Let us first remark that $L(n,-n+2) =L(n,-n+1) =L(n,-n) = \ldots$, so we will focus on the values of $\theta$ from $n$ to $-n+2$. An example of all $L(n,\theta)$ for  $n=3$  is given in Figure \ref{fig:04}.
\begin{figure}[h]
$${\epsfbox{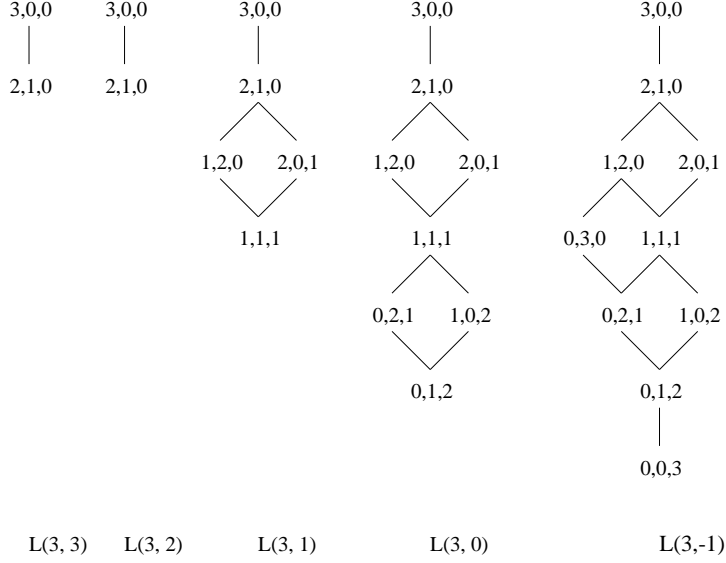}}$$
\caption{The lattices $L(3,\theta)$ with $\theta \in [-1,3]$}
\label{fig:04}
\end{figure}
\subsection{Lattice structure of $L(n,\theta)$ and characterization of its elements}
The purpose of this section is to show that each $L(n,\theta)$ is a lattice and that the set of all $L(n,\theta)$, $n\geq \theta \geq -n+2$, forms an increasing sequence of lattices  from the lattice $L(n)$ to the lattice $L(-n+2)$, this last one containing all compositions of $n$. Moreover we also show that $L(n,2) = SPM(n) \unlhd L_B (n) \unlhd L(n,1)$ where $\unlhd$ denotes the suborder relation. Before proving that each $L(n,\theta)$ is a lattice, let us first show that it is an order and study this order by means  of the energy of each composition. Let us recall that the energy $E$ of a composition $a$ is defined by $E(a) = \sum_{i=1}^{i=n} (i-1) a_i$ \cite{GK93}. It is clear to see that $E(b) = E(a) +1$ if $b$ is obtained from $a$ by applying a transition $T(\theta)$, so the set $L(n,\theta)$ induces an order, where $a \geq_{\theta} b$ if we can apply a sequence of transitions $T(\theta)$ from $a$ to obtain $b$. When $\theta \geq 2$, the obtained compositions are all decreasing ones, so we will only consider the strictly positive values of the composition and not the sequence of $0$s ending. So from right now, without indication, in the case $\theta \geq 2$, we will call by $a$ the first part of strictly positive integers of a given composition $a$.

Let us first introduce a new notion of sub-sequence for a composition: we say that a sequence of length $l+2$, $l \geq 1$, in $L(n,\theta)$ is a $(l,\theta)$-sequence if it has the following form:
$$k(k-\theta +2)(k-2\theta +3)(k-3\theta +4)\ldots (k- l\theta +l+1)(k-(l+1)\theta +l+3)$$
We can now give the characterization of $L(n,\theta)$ 
\begin{theo} 
The necessary and sufficient conditions $C(\theta)$ for $a$ to belong to $L(n,\theta)$ are:

i) $\forall i, 1\leq i \leq n-1,$ $a_i - a_{i+1} \geq \theta -2$ 

ii) $a$ does not contain any $(l,\theta)$-sequence.
\end{theo}
\Proof

{\em Necessary condition:} 
 It is easy to see that $N$ satisfies condition $C(\theta)$. Assume that $ a \in L(n,\theta)$, condition i) comes directly from the definition. Assume now that $a$ contains a $(l,\theta)$-sequence. The case $l=1$ can not occur, since in this case, it would not been possible to move back any grain from the piles of the sequence $(k, k-\theta+2, k- 2\theta +4)$ and so go back to $N$, and then $a$ would not be reachable from $N$. We can then suppose that $l \geq 2$. Our intention is to find another composition of $L(n,\theta)$ containing a $(l',\theta)$-sequence with $l'<l$. By repeatedly applying this construction, we will obtain a contradiction. Since $a \in L(n,\theta)$, there exists a path of inverse transitions (respecting $T(\theta)$) from $a$ to $N$. Let us denote by $a'$ the first composition on this path from $a$ to $N$ obtained by modification of the $(l,\theta)$-sequence of $a$. $a' \in L(n,\theta)$ so $a'$ satisfies i) and then this modification is neither a decrease of the first pile nor a growth of the last pile, so $a'$ contains a $(l',\theta)$-sequence with $l' <l$, which was what we wanted to show. 

{\em Sufficient condition}: Let $a$ be a composition satisfying condition $C(\theta)$. Our intention is to find a path of inverse transitions of $T(\theta)$ from $a$ to $N$. To do this, it is sufficient to show that there exists a composition $a'$ satisfying $C(\theta)$ such that $a' \geq_{\theta} a$. If $\forall i$ $ a_i - a_{i+1} \geq \theta -1$, then $a'$ obtained from $a$ by applying the inverse transition of $T(\theta)$ at the position 1 will satisfy $C(\theta)$. Otherwise, let $i$ be the first index such that   $a_i - a_{i+1} = \theta -2$. Since $a$ satisfies condition $C(\theta)$ then $a_{i-1} - a_i \geq \theta -1$ and $a_{i+1} - a_{i+2} \geq \theta -1$. Let us define the composition $a'$ as follows:
$$\forall j \quad  a'_j  = \left\{\begin{array}{lll}
			   a_i + 1  & \mbox{ if } j = i   \\
			   a_{i+1} -1 & \mbox{ if } j=i +1 \\
			   a_j & \mbox{otherwise.}
			\end{array}
		           \right. 
$$   
This composition $a'$ clearly satisfies condition $C(\theta)$, and we have then proved the theorem. 
\hfill$\Diamond$\\   
The following corollary is immediate:
\begin{coro}
Let $a$ be an element of $L(n,\theta)$, for all integers $i,l$ such that $i+l \leq n$, we have $a_i - a_{i+l} >l(\theta-1)-2$. 
\end{coro}

 By using the previous result about the characterization of elements of $L(n,\theta)$, we can now study the nature of the order $L(n,\theta)$.
\begin{theo}
Let $a$ and $b$ be two elements of $L(n,\theta)$, then $a \geq_{\theta} b$ if and only if $a \geq_{D} b$ by dominance ordering extended to all compositions of $n$, i.e. $\forall i, 1\leq i\leq n$, $\sum_{j=1}^{j=i} a_j \geq  \sum_{j=1}^{j=i} b_j$.
\end{theo}
\Proof To prove that $a \geq_{\theta} b \Rightarrow a \geq_{D} b$, let us first consider a transition $T(\theta)$ $x \rightarrow y$ at the position $i_0$ on two elements on a sequence from $a$ to $b$. It is easy to see that:
$$\forall i \quad \sum_{j=1}^{j=i} y_j   = \left\{\begin{array}{ll}
			\sum_{j=1}^{j=i} x_j - 1 & \mbox{ if } i = i_0   \\
			\sum_{j=1}^{j=i} x_j    & \mbox{otherwise}
			\end{array}
		 \right. 
$$
which shows that $x \geq_{D} y$, so $a \geq_{D} b$ by transitivity.

On the other hand, let $a, b \in L(n,\theta)$ and $a \geq_{D} b$, our purpose is to construct a sequence of transitions $T(\theta)$ from $a$ to $b$. We will first find a composition $a'$ such that 
$a\geq_{\theta} a' \geq_{D} b$. Let us consider the first index $i$ such that $a_i > b_i$. As $a \geq_D b$, we have $a_j = b_j$ $\forall j<i$. Let $l$ be the smallest integer such that $\sum_{j=1}^{j=i+l} a_j =\sum_{j=1}^{j=i+l} b_j$. We have then  $\sum_{j=1}^{j=i+l-1} a_j > \sum_{j=1}^{j=i+l-1} b_j$, so $a_{i+l} < b_{i+l}$. If we can apply a transition $T(\theta)$ on a position between $i$ and $i+l-1$ of $a$, the obtained composition $a'$ will satisfy the condition $a\geq_{\theta} a' \geq_{D} b$. On the other hand, if a transition $T(\theta)$ can not be applied, then $\forall j, i\leq j \leq i+l-1 $ $  a_j -a_{j+1} < \theta$, which implies $a_i - a_{i+l} \leq l(\theta -1)$, and $b_i -b_{i+l} \leq l(\theta-1) -2$, so $b \not \in L(n,\theta)$ by Corollary 1, which is a contradiction.
We have then found the composition $a'$ and, by continuing this construction, we obtain $a \geq_{\theta} b$. 
\hfill$\Diamond$\\
We can now prove that $L(n,\theta)$ is a lattice:
\begin{theo}
$L(n,\theta)$ is a lattice, where
$$inf_{\theta}(a,b)= c \mbox{ if } \mbox{ if } c_i = min (\sum_{j=1}^{j=i} a_j, \sum_{j=1}^{j=i} b_j) - \sum_{j=1}^{j=i-1} c_j.$$
\end{theo}
\Proof
Let $a$ and $b$ be two elements of $L(n,\theta)$, and let $c$ be the composition which is defined as follows:
$$c_i = min (\sum_{j=1}^{j=i} a_j, \sum_{j=1}^{j=i} b_j) - \sum_{j=1}^{j=i-1} c_j.$$
Let us prove that $c$ satisfies condition $C(\theta)$. Suppose first  that $c$ contains a $(l,\theta)$-sequence, i.e.:
$$c_i, \ldots, c_{i+l+1} = k(k-\theta +2)(k-2\theta +3)(k-3\theta +4)\ldots (k- l\theta +l+1)(k-(l+1)\theta +l+3)$$ 
 Knowing that $\forall m, 1\leq m\leq n$, $\sum_{j=1}^{j=m} c_j = min (\sum_{j=1}^{j=m} a_j, \sum_{j=1}^{j=m} b_j)$, without loss of generality we may suppose that $\sum_{j=1}^{j=i-1} c_j= \sum_{j=1}^{j=i-1} a_j$. Two cases are now possible: either $\sum_{j=1}^{j=i} c_j= \sum_{j=1}^{j=i} a_j$, or $\sum_{j=1}^{j=i} c_j= \sum_{j=1}^{j=i} b_j$. In the first case, we obtain $a_i = c_i$, and as  $a$ satisfies condition $C(\theta)$, $\sum_{j=1}^{j=i+l+1}a_j < \sum_{j=1}^{j=i+l+1} c_j$, which is a contradiction. The other case is analogous. This shows that $c$ satisfies condition ii) of Theorem 6. Condition i) can be proved in the same way. So $c$ belongs to $L(n,\theta)$ and is equal to $inf_{\theta}(a,b)$. Moreover, $L(n,\theta)$, being an order with the maximal element $N$, is then  a lattice.
\hfill$\Diamond$\\
Furthermore,  if we consider all $L(n,\theta)$ with $n \geq \theta \geq -n+2$, we obtain an interesting result on the relation between them:
\begin{theo}
The lattices $L(n,\theta)$ form an increasing sequence with regard to the suborder relation:
$$ L(n,n)\unlhd \ldots \unlhd L(n,2) \unlhd L(n,1) \unlhd L(n,0) \unlhd \ldots \unlhd L(n,-n+2)$$
Moreover,   $|L(n,-n+2)| = \left( \begin{array}{c} 2n-1\\n \end{array} \right)$
\end{theo}
\Proof
The first part of the theorem immediately comes from the fact that $\forall \theta, -n+3 \leq \theta \leq n$,  $a \leq_\theta b \Leftrightarrow  a\leq_D b \Leftrightarrow a\leq_{\theta -1} b$. Let us now consider the lattice $L(n,-n+2)$, it is easy to see that every composition satisfies condition $C(-n+2)$, so this lattice contains all compositions and its cardinal is then equal to $\left( \begin{array}{c} 2n-1\\n \end{array} \right)$.
\hfill$\Diamond$\\
Let us complete this sequence by the following immediate result, where $\unlhd$ still denotes the suborder relation.
\begin{propo}
$L(n,2) =SPM(n) \unlhd L_B (n) \unlhd L(n,1)$.
\end{propo}
\subsection{Fixed point and length of maximal chains of $L(n,\theta)$}
 In \cite{GK93}, while studying $SPM(n)$, Goles and Kiwi \cite{GK93} have presented a formula characterizing its fixed point. Here, using Theorem 5, we will also find the formula for the fixed point of each $L(n,\theta)$. Note that since $L(n,\theta)$ is a lattice, it can only have one fixed point which is the smallest element of the lattice. Let us first remark that for a given $\theta$, any integer $n$ can be uniquely written as follow:
$$n= \frac{k(k+1)}{2} |1-\theta| + l(k+1) +p $$
where $ 0 \leq l < |1-\theta|$, $0 \leq p \leq k+1$.\\
For $\theta \leq 1$, let
$$P(\theta) =  0,\ldots,0,l,l+(1-\theta), l+2(1-\theta),\ldots,l+(k-p)(1-\theta),$$
$$ l+(k-p+1)(1-\theta)+1,\ldots,l+k(1-\theta)+1. $$
For $\theta \geq 2$, let 
$$P(\theta) = l+k(\theta-1), l+(k-1)(\theta-1),\ldots,l+p(\theta-1),$$
$$l+(p-1)(\theta-1)+1,\ldots,l+(\theta-1)+1,l+1.$$
 It is clear that $P(\theta)$ satisfies the condition $C(\theta)$, so $P(\theta) \in L(n,\theta)$. It is also clear to see that we cannot apply the transition $T(\theta)$ over $P(\theta)$, thus we have the following result:
\begin{propo}
$P(\theta)$ is the fixed point of $L(n,\theta)$ 
\end{propo}
Here,  we can find the formula for the fixed point in the particular case  of $SPM(n)$, i.e.
$$P(2) = k,k-1,\ldots,p+1,p,p,p-1,\ldots,2,1$$
which was given in \cite{GK93}.

On the other hand, the problem of the length of maximal chains is also worth further inverstigation. In $L(n,\theta)$, $a\rightarrow b \Rightarrow E(b) = E(a) +1$, so all chains between two elements $a$ and $b$ of $L(n,\theta)$ have the same length which is equal to $C(a,b) = E(b) -E(a)$. So, in $L(n,\theta)$, any maximal chains' length is equal to the length of any chain from $N$ to $P(\theta)$, and then
\begin{propo}
The maximal chain lenght in $L(n,\theta)$ for $n= \frac{k(k+1)}{2} |1-\theta| + l(k+1) +p $ is equal to $E(P(\theta))$, that is\\
If $\theta \geq 2$:
$$(\theta-1)\frac{(k-1)k(k+1)}{6} + l\frac{k(k+1)}{2} + p \frac{2k-p+1}{2},$$
and if $\theta \leq 1$:
$$(1-\theta)\frac{(3n-k-2)k(k+1)}{6} + l\frac{(k+1)(2n-k-2)}{2} + p\frac{2n-p-1}{2}.$$
\end{propo}
\Proof
Let us consider the case $\theta \geq 2$, the other case being proved similarly.
Since $$P(\theta) = l+k(\theta-1), l+(k-1)(\theta-1),\ldots,l+p(\theta-1),$$
$$l+(p-1)(\theta-1)+1,\ldots,l+(\theta-1)+1,l+1,0,\ldots,0$$
then $$E(P(\theta)) = \sum_{i=1}^{i=n} (i-1) P_i$$
$$= \sum_{i=1}^{i=k}(i-1)(k-i+1)(\theta-1) + \sum_{i=1}^{i=k+1}(i-1)l + \sum_{i=k-p+2}^{i=k+1} (i-1)$$
$$= (\theta-1)\frac{(k-1)k(k+1)}{6} + l\frac{k(k+1)}{2} + p \frac{2k-p+1}{2}.$$
\hfill$\Diamond$\\

\bibliographystyle{plain}
\bibliography{Spiles}

\end{document}